\magnification=1000
\baselineskip=15pt
\hsize=125truemm
\vsize=190truemm
\hoffset=4truemm
\voffset=5truemm
\hfuzz=10pt
\headline={\hfil}
\footline={\hfil\tenrm\folio\hfil}
\font\obf=cmbx10 scaled\magstep1

\parindent=5mm
\centerline{\sevenrm Journal of Pure and Applied Algebra 171 (2-3) (2002) 171-184}
\bigskip

\centerline{\obf On The Class Group of a Graded Domain }
\vskip1.5cm
 
\centerline{ S. El Baghdadi$^{a}$, L. Izelgue$^{b}$, S. Kabbaj$^{c,}$}
\bigskip

\centerline {\sevenrm $^{a}$ Department of Mathematics, FST, P.O. Box 523,
Beni Mellal, Morocco.}

\centerline {\sevenrm $^{b}$ Department of Mathematics, University of Marrakech, P.O. Box S15,
Marrakech, Morocco. }

\centerline {\sevenrm $^{c}$ Department of Mathematical Sciences, KFUPM, P.O.Box 5046,
Dhahran 31261, Saudi Arabia.}

\vskip1cm\hrule\medskip
\noindent{\bf Abstract}
\medskip
This paper studies the class group of a graded integral 
domain $R=\oplus_{\alpha \in \Gamma}R_{\alpha}$. We prove that if the extension 
$R_0 \subset R$ is inert, then $Cl(R)=HCl(R)$ if and only if $R$ is almost normal. 
As an application, we state a decomposition theorem for class 
groups of semigroup rings, namely, $Cl(A[\Gamma])\cong Cl(A)\oplus HCl(K[\Gamma])$ if 
and only if $A[\Gamma]$ is integrally closed. This recovers the well-known results developed 
for the classic contexts of polynomial rings and Krull semigroup rings. Further, we obtain
an interesting result on the natural homomorphism $\phi: Cl(A)\rightarrow Cl(A[\Gamma])$, that is, 
$Cl(A[\Gamma])=Cl(A)$ if and only if $A$ and $\Gamma$ are integrally closed and $Cl(\Gamma)=0$. 
Our results are backed by original examples.


\medskip\hrule\vskip1.5cm

\noindent {\bf 0. Introduction}
\medskip

All rings considered in this paper are integral domains. 
  Throughout, $\Gamma$ will always denote a torsionless 
  grading monoid. That is, $\Gamma$ is a commutative cancellative 
monoid, written additively, and the quotient group generated by $\Gamma$ is a
 torsion-free abelian group. By a graded domain  
$R=\oplus_{\alpha \in \Gamma}R_{\alpha}$, we mean an integral domain $R$ graded by 
an arbitrary torsionless grading monoid $\Gamma$. Suitable background on 
torsionless grading monoids and $\Gamma $-graded rings is [21].\par

 Let $R$ be an integral domain. Following [6], we define the class group of 
$R$, denoted $Cl(R)$, to be the group of $t$-invertible fractional $t$-ideals of 
$R$ under $t$-multiplication modulo its subgroup of principal fractional ideals. 
Divisibility properties of a domain $R$ are often reflected in group-theoretic 
properties of $Cl(R)$. For $R$ a Krull domain, $Cl(R)$ is the usual divisor 
class group of $R$. In this case, $Cl(R)=0$ if and only if $R$ is factorial. 
If $R$ is a Pr\"ufer domain, then $Cl(R)=Pic(R)$ is the ideal class group of $R$. 
In this case, $Cl(R)=0$ if and only if $R$ is a Bezout domain. We assume familiarity 
with class groups and related concepts, as in [6], [10], and [12].\par  

On one hand, a well-known result is that if $R$ is a $Z_+$-graded Krull domain, then $Cl(R)$ 
is generated by the classes of homogeneous height-one prime ideals of $R$ [10, Proposition 10.2], 
i.e., $Cl(R) = HCl(R)$, where $HCl(R)$ is the homogeneous  class group of $R$. 
In [3, Theorem 4.2], D.F. Anderson showed that the same holds for any $\Gamma$-graded Krull 
domain, where $\Gamma$ is an arbitrary torsionless grading monoid. So one may remove 
the ``Krull assumption" and legitimately ask the following: {\it For an arbitrary graded domain $R$, how 
does the equality ``$Cl(R)= HCl(R)$" reflect in ring-theoretic properties of $R$?}\par

On the other hand, many authors investigated the problem of characterizing ring-theoretic properties 
in terms of Picard groups. In [7], Bass and Murthy  proved that for an integral domain $A$, if 
$Pic(A[X, X^{-1}])=Pic(A)$ then $A$ is seminormal. However, Pedrini showed that 
the converse fails to be true in general [22, p. 96]. Many years later, Gilmer and 
Heitmann [14] solved completely the problem of characterizing ``seminormality" in terms of Picard groups.
They stated that $Pic(A[X])=Pic(A)$ if and only if $A$ is seminormal. In the same line, in 1982, 
the Andersons [1] examined the property of almost normality for graded domains. 
They established that if $R$ is an almost normal graded domain with $R_0\subset R$ inert, 
then $Pic(R)=HPic(R)$. So the problem remained somehow open. However, in 1987, 
Gabelli proved that for an integral domain $A$, $Cl(A[X])=Cl(A)$ if and only if $A$ is 
integrally closed [11, Theorem 3.6]. Recall for convenience that $A[X]$, graded in the natural way, 
is almost normal if and only if $A[X]$ (and hence $A$) is integrally
closed. This motivates our second question: {\it For an arbitrary graded domain $R$, how 
does ``almost normality" reflect in group-theoretic properties of $Cl(R)$?} \par 

This paper contributes to the study of class groups of graded integral domains. 
It particularly provides a satisfactory (and 
unique) answer to the previous two questions. As an application, we state a decomposition 
theorem for class groups of semigroup rings. Indeed, Section 1 examines the interconnection 
between ``almost normality" and the equality ``$Cl(R)=HCl(R)$" for a graded domain $R$. 
More precisely, we show, in Theorem 1.1, that if $R_0 \subset R$ is inert, then $Cl(R)=HCl(R)$ if and only if $R$ 
is almost normal. Some interesting contexts for this result are $Z_+$-graded domains and polynomial rings.
However, Example 1.11 illustrates its failure if one omits the condition  ``$R_0 \subset R$  inert". 
In the first part of Section 2, we focus on the specific case of semigroup rings, which provide an important 
class of graded domains. We establish the following decomposition theorem, Theorem 2.7: 
For an integral domain $A$, with quotient field $K$,  $Cl(A[\Gamma])\cong Cl(A)\oplus HCl(K[\Gamma])$ if and only if 
$A[\Gamma]$ is integrally closed. This recovers most of the previous results stated 
for the classic contexts of polynomial rings [11] and Krull semigroup rings [3]. The second part 
of Section 2 is devoted to semigroups. Here we extend Chouinard's 
results on Krull semigroups (cf. [8]) to arbitrary semigroups. As an application, we establish
an interesting result, Theorem 2.12, on the natural homomorphism 
$\phi: Cl(A)\rightarrow Cl(A[\Gamma])$, that is, $Cl(A[\Gamma])=Cl(A)$ if
and only if $A$ and $\Gamma$ are integrally closed and $Cl(\Gamma)=0$. 
\vskip1cm 

\noindent {\bf 1. The class group of a graded domain}
\medskip  

The discussion which follows, concerning basic facts and notations connected with graded domains, 
 will provide some background 
 to the main result of this section and will be of use in its proof. In this section, 
 $ R=\oplus_{\alpha \in \Gamma}R_{\alpha}$
denotes a $\Gamma$-graded domain and $S$ the multiplicatively closed subset of all 
nonzero homogeneous elements of $R$.
Thus $R_S$ is a  $<\Gamma>$-graded domain with
$$(R_S)_{\alpha}=\{{a\over b}; \, a\in R_{\beta}, \,0\not=b \in R_{\gamma},\hbox{and}\,
 \alpha=\beta-\gamma\}.$$ 
 In particular, $(R_S)_{0}$ is a field and each nonzero
homogeneous element of $R_S$ is a unit. It is well-known that the domain $R_S$, often called the 
 homogeneous quotient field of $R$,  is a completely integrally closed 
GCD-domain [2, Proposition 2.1]. If $R$ is $Z_{-}$ or $Z_{+}$-graded, then
 $R_S\cong (R_S)_{0}[X,X^{-1}]\cong (R_S)_{0}[Z]$. Semigroup rings $A[\Gamma]$, graded in the natural 
way with $degX^{\alpha}=\alpha$, constitute perhaps the most important class 
of  $\Gamma$-graded domains. The homogeneous quotient field of $A[\Gamma]$
  is the group ring $K[G]$, where $K$ is the quotient field of $A$ and $G=<\Gamma>$, the 
  quotient group of $\Gamma$. \par 

We say that the graded domain $R$ is almost normal [1] if for each homogeneous element 
$x\in R_S$ of nonzero degree which is integral over $R$ is
actually in $R$. Any integrally closed graded domain is almost normal; moreover, it is 
well-known that $R$ is integrally closed if and only if $R$ is almost normal and $R_0$
is integrally closed in $(R_S)_0$. On the other hand, following Cohn [9], we say that an
extension $A\subset B$ of integral domains is inert if whenever $xy\in A$ for some $x,
y\in B$, then $x=ru$ and $y=su^{-1}$ for some $r,s\in A$ and $u$ a unit of $B$. The 
extension $R_0\subset  R=\oplus_{\alpha \in \Gamma}R_{\alpha}$ is inert if and only 
if for each $\alpha \in \Gamma\cap
(-\Gamma)$, $R_{\alpha}$ contains a unit. This, particularly, happens if either 
$R=R_0[\Gamma]$ is a semigroup ring, or $R_0$ is a field, or $R$ is $\Gamma$-graded 
with $\Gamma\cap (-\Gamma)=0$.\par

Now, Let's review  some terminology related to the $v-$ and $t-$operations . Let $A$ be 
any domain with quotient field $K$. By 
an ideal of $A$ we mean an integral ideal of $A$. Let  $I$ and $J$ be two  nonzero 
fractional ideals of $A$. We define the fractional ideal $(I:J)=\{x\in K\mid xJ\subset I\}$.
 We denote $(A:I)$ by $I^{-1}$ and  $(I^{-1})^{-1}$ by $I_v$.
We say that $I$ is divisorial  or  a $v$-ideal of $A$ if $I_v=I$. The ideal $I$ is $v$-finite if 
$I=J_v$ for some finitely generated fractional ideal $J$ of $A$.
 For a nonzero fractional ideal $I$ of $A$, we define $I_t=\cup\{J_v\mid\, J\subset I$ 
finitely generated$\}$. The ideal $I$ is a $t$-ideal if $I_t=I$. Under the operation 
$(I,J)\mapsto (IJ)_t$, the set of $t$-ideals of $A$ is a semigroup with unit $A$. An   
invertible element for this operation is called a $t$-invertible $t$-ideal of $A$.
  For more details about these notions, see [12, Sections 32 and 34].\par 

A fractional ideal $I$ of  the graded domain $R$ is homogeneous if there exists a
 nonzero homogeneous element $s$ of $R$ such that $sI$ is a homogeneous (integral) ideal of $R$.
 Each homogeneous fractional ideal of $R$ is contained in $R_S$. Moreover, if $I$ and 
$J$ are  nonzero homogeneous fractional  ideals of $R$, then $(I:J)$ is also a 
homogeneous fractional ideal of $R$, and so are $I^{-1}$ and $I_v$ [2, Proposition 2.5]. 
Let $T(R)$ (resp., $HT(R)$) denote the  group of   
$t$-invertible fractional $t$-ideals (resp., homogeneous  $t$-invertible fractional $t$-ideals) 
of $R$, and $P(R)$ (resp., $HP(R)$), its subgroup of principal fractional ideals. 
We have  $Cl(R)=T(R)/P(R)$ and  $HCl(R)=HT(R)/HP(R)$,  a subgroup of $Cl(R)$. Now, 
let $x\in R_S$ with $x=x_{\alpha_1}+\cdots+x_{\alpha_n}$, where
$x_{\alpha_i}\in (R_S)_{\alpha_i}$ and $\alpha_1<\cdots<\alpha_n$. We define the 
content of $x$, denoted  $C(x)$, to be $C(x)=(x_{\alpha_1},\ldots,x_{\alpha_n})$, 
the homogeneous fractional ideal of $R$ generated by the homogeneous components of $x$ 
 in $R_S$. If  $I\subset R_S$ is a fractional ideal of $R$, then $I$ is homogeneous
if and only if  $C(x)\subset I$ for each $x\in I$.
A well-known result due to Northcott   [20] states  that, for each $x,y\in R$, 
 $C(x)^nC(xy)= C(x)^{n+1}C(y)$ for some integer
$n\ge 0$.\par 
\vskip1cm

We now announce our main result of this section. It sheds light on the interconnection 
between ``almost normality" and the equality ``$Cl(R)=HCl(R)$" for a graded domain $R$. 
\bigskip

\noindent{\bf Theorem 1.1.} \quad {\sl Let $R=\oplus_{\alpha \in \Gamma}R_{\alpha}$ be  a
$\Gamma$-graded domain such that $R_0\subset R$ is inert. Then $Cl(R) = HCl(R)$ if
and only if $R$ is almost normal.}\medskip

The proof of this theorem requires some preliminaries.
\bigskip

Let $I$ be a fractional ideal of $R$ and assume that there exists $s\in S$ such that
$sI\subset R$. We define the content of $I$, denoted $C(I)$, to be the homogeneous
fractional ideal of $R$ generated by the homogeneous components in  $R_S$ of all elements
of $I$ . We have $C(I)= \sum_{x\in I}C(x)$, and $I$ is homogeneous if and only if $C(I)=I$. 
The next two lemmas deal with technical properties of the content of a fractional ideal. 
\bigskip

 \noindent{\bf Lemma 1.2.}\quad {\sl Let $I$ be a fractional ideal of $R$ with $sI\subset R$ 
 for some $s\in S$. Then\par
(1) $C(HI)=HC(I)$, for each homogeneous fractional ideal $H$ of $R$.  \par
(2) $C(I_v)_v=C(I)_v$. \par
(3) $C(I_t)_t=C(I)_t$.}
\medskip 
\noindent {\bf Proof.}\quad (1) is straightforward.\par 
(2) We have $I\subset C(I)$, hence $I_v\subset C(I)_v$ and $C(I_v)_v\subset C(I)_v$.
The reverse inclusion is trivial. \par
 (3) We first show that $I_t$ is homogeneous whenever $I$ is. For,  let $x\in I_t$. Then there exists a finitely generated fractional ideal $F\subset I$ such that $x\in F_v$.
Since $I$ is homogeneous, $C(F)\subset I$. Hence by (2),  $C(x)\subset C(F_v)\subset C(F)_v\subset I_t$.
Therefore, $I_t$ is homogeneous. Now, since $I_t\subset C(I)_t$ and  $C(I)_t$ is  
homogeneous, then $C(I_t)\subset C(I)_t$. Hence $C(I_t)_t\subset C(I)_t$.  
The reverse inclusion is trivial. $\diamondsuit$ 
\bigskip

\noindent{\bf Lemma 1.3.}\quad {\sl Let $x_1,\ldots, x_n\in R$ such that $(C(x_1)+\cdots+C(x_n))_v=R$, and $a\in S$.\par 
Then  $(a, x_1,\ldots, x_n)_v=R$.}
\medskip

\noindent {\bf Proof.}\quad Let $u\in qf(R)$ such that $a, x_1,\ldots, x_n\in uR$. Then
$u={a\over r}$ for some $r\in R$ and, for each $i$, $x_i={a\over r}r_i$ for some $r_i\in
R$.  By  the Northcott's result, for each $i=1,\ldots,n$, there exists a positive
integer $N_i$ such that $C(x_i)^{N_i}C(rx_i)=C(x_i)^{N_i+1}C(r)$. Let
$N=N_1+\cdots+N_n$. Then $C(x_i)^{N}C(rx_i)=C(x_i)^{N+1}C(r)$, for each $i=1,\ldots,n$. On the other hand, for each $i=1,\ldots,n$, $C(rx_i)=C(ar_i)=aC(r_i)\subset aR$, hence $C(x_i)^{N+1}C(r)\subset aR$. Therefore
$({\displaystyle\sum_{i=1}^n}C(x_i)^{N+1})C(r)\subset aR$. Since $({\displaystyle\sum_{i=1}^n}C(x_i))_v=R$, then
$({\displaystyle\sum_{i=1}^n}C(x_i)^{N+1})_v=R$. It follows that $C(r)_v\subset aR$ and
hence ${r\over a}\in R$. Thus $1\in uR$ and hence $(a, x_1,\ldots, x_n)_v=R$, as desired.  $\diamondsuit$ 
\bigskip 

Next we give a  criterion for $t$-invertibility in a graded
domain (see also [15, section 4]).\bigskip

 \noindent{\bf Lemma 1.4.}\quad {\sl Let $I$
be an ideal of $R$ such that $C(I)$ is $t$-invertible. Then $I$ is $t$-invertible if
 and only if $I_tR_S$ is principal.}
\medskip

 \noindent {\bf Proof.}\quad If $I$ is  $t$-invertible, by [5, Proposition 2.2], $I_tR_S$ is a
$t$-invertible $t$-ideal. Since $R_S$ is a GCD-domain, $I_tR_S$ is principal. Conversely, set
$J=C(I)^{-1}I_t$. By Lemma 1.2, $C(J)_t=C(C(I)^{-1}I)_t=(C(I)^{-1}C(I))_t=R$. Further, since
$JR_S=I_tR_S$ is principal, it suffices to show  that if $C(I)_t = R$ and
$IR_S$ is principal, then $I$ is $t$-invertible. Let $x_1, \ldots, x_n\in I$ such that
$({\displaystyle\sum_{i=1}^n}C(x_i))_v=R$.  Since $IR_S$ is principal, then $IR_S=aR_S$ for 
some  $a\in I$. Now,  if $x\in I$, then $x={{ar}\over s}$
for some $r\in R$ and $s\in S$. Thus $x(s,x_1,\ldots, x_n)_v\subset (a,x_1,\ldots,
x_n)_v$. By Lemma 1.3, $(s,x_1,\ldots, x_n)_v=R$, so $x\in (a,x_1,\ldots, x_n)_v$ 
and hence $I\subset (a,x_1,\ldots, x_n)_v$. On the other hand,  there
exists $t\in S$ such that for each $i$, $x_i={{ar_i}\over t}$ for some  $r_i\in R$.
Hence $tI\subset(ta,tx_1,\ldots,tx_n)_v\subset aR$, i.e., ${t\over a}\in I^{-1}$.
Therefore, $t=a{t\over a}\in II^{-1}$ and, by Lemma 1.3, 
$R=(t, x_1,\ldots, x_n)_v\subset (II^{-1})_t$. Hence, $I$ is $t$-invertible.   $\diamondsuit$ \bigskip

Notice that a useful case of Lemma 1.4 is when $C(I) =R$.\bigskip

\noindent{\bf Lemma 1.5.}\quad {\sl Let $a\in R_S$ be a nonzero element and $P_a=aR_S\cap R$.
Then $P_a=uJ$ for some $u\in R_S$ and  some homogeneous ideal $J$ of $R$ if and only if $P_a=aC(a)^{-1}$.}
\medskip

\noindent {\bf Proof.}\quad If $P_a=uJ$, then $P_aR_S=aR_S=uR_S$; so
there exist $s,t\in S$ such that $u={s\over t}a$. Now, since ${s\over t}aJ\subset R$, then
${s\over t}J\subset C(a)^{-1}$, and hence $P_a={s\over t}aJ\subset aC(a)^{-1}$. The
reverse inclusion is trivial. Conversely, if $P_a=aC(a)^{-1}$, let $s\in C(a)$ be a 
nonzero homogeneous element. Then $P_a={a\over s}(sC(a)^{-1})$, and hence we may take
$J=sC(a)^{-1}$. $\diamondsuit$ 
\bigskip

We next state our key lemma. Its main effect is to link,
under the stated hypothesis, $t$-invertibility to almost normality.\bigskip

\noindent{\bf Lemma 1.6.}\quad  {\sl Assume that $R_0\subset R$ is inert. The following
statements are equivalent.\par
 (i) $R$ is almost normal;\par 
(ii) For each $v$-finite $v$-ideal $I$ of $R$, there exist $u\in R_S$ and a homogeneous
$v$-finite $v$-ideal $J$ of $R$ such that $I=uJ$;\par
 (iii)  For each $t$-invertible $t$-ideal $I$ of $R$, there exist $u\in R_S$ and a
homogeneous $t$-invertible  $t$-ideal $J$ of $R$ such that $I=uJ$.}\medskip   

\noindent {\bf Proof.}\quad The equivalence (i)$\Leftrightarrow $(ii) follows from [1, Theorem 3.7(2) and Theorem 3.2].\par 
(ii)$\Rightarrow$ (iii) is obvious.\par
(iii)$\Rightarrow$ (i) Let $a\in R_S$ be  homogeneous of nonzero degree and integral
over $R$. Set $P_a=(1-a)R_S\cap R$ and let $f(X)=X^n+r_{n-1}X^{n-1}+\cdots+r_0\in
R[X]$ such that $f(a)=a^n+r_{n-1}a^{n-1}+\cdots+r_0=0$. Regrouping terms of the same
degree, we may assume that the $r_i$'s are homogeneous of pairwise distinct nonzero
degrees. We have $f(X)=(X-a)g(X)$, where  $g(X)=X^{n-1}+b_ {n-2}X^{n-2}+\cdots+b_0$
with the $b_i$'s are homogeneous elements of $R_S$ of pairwise distinct nonzero
degrees. On the other hand, $f(1)=(1-a)g(1)\in P_a$, so  $1+r_{n-1}+\cdots+r_0\in
P_a$, moreover $C(1+r_{n-1}+\cdots+r_0)=R$. It follows that $C(P_a)=R$. Since $P_aR_S=(1-a)R_S$ 
is principal and $P_a$ is a $t$-ideal,  then $P_a$ is
$t$-invertible (cf. Lemma 1.4).  Therefore, there exist $u\in R_S$ and $J$ a homogeneous 
$t$-invertible $t$-ideal such that $P_a=uJ$. By Lemma 1.5, $P_a=(1-a)C(1-a)^{-1}$. We deduce 
that $f(1)=(1-a)(1+b_{n-2}+\cdots+b_0)\in
(1-a)C(1-a)^{-1}$, i.e., $1+b_{n-2}+\cdots+b_0\in (1,a)^{-1}$. It follows that 
$a+ab_{n-2}+\cdots+ab_0\in R$, hence $a\in R$ since the  $b_i$'s are homogeneous of
pairewise distinct nonzero degrees. $\diamondsuit$ \bigskip 

{\bf Proof of Theorem 1.1.} It follows from Lemma 1.6.  $\diamondsuit$ 
\bigskip

\noindent{\bf Remark 1.7.}\quad (1) In [17], Matsuda constructed an example illustrating 
the fact that in [1, Theorem 3.7 (2)] and hence in Lemma 1.6 (i)$\Rightarrow$ (ii)	the ``$R_0\subset R$  inert"
hypothesis cannot be deleted.\par
(2) In  Lemma 1.6, we need this hypothesis   only for the 
 implication (i)$\Rightarrow$ (ii). However, [1, Theorem 3.2 and Theorem 3.7 (1)]
shows that we can omit it  if we assume that $R$ contains a (homogeneous)
unit of nonzero degree. In this case, $R$ is almost normal if and only if $R$ is
integrally closed.\par
(3) By substituting the hypothesis ``$R$ contains a   unit of nonzero
degree"  for ``$R_0\subset R$ inert", the statement of Theorem 1.1 remains true,
that is, $Cl(R)=HCl(R)$ if and only if $R$ is integrally closed.\bigskip   

\noindent {\bf Corollary 1.8.} \quad {\sl If $R=R_0 \oplus R_1\oplus ...$ is $Z_+$-graded,
then $Cl(R) = HCl(R)$ if and only if $R$ is almost normal.}  $\diamondsuit$ 
\bigskip

Let $A$ be an integral domain and $X$ an
indeterminate over $A$. Using [1, Proposition 5.8], one may show that $HCl(A[X])=Cl(A)$. 
Thus we reobtain Gabelli's result:\bigskip

\noindent {\bf Corollary 1.9.} [11, Theorem 3.6] \quad  {\sl $Cl(A[X])=HCl(A[X])(=Cl(A))$ if and only if $A$ is integrally
closed.} $\diamondsuit$ \bigskip

 Let $A\subset B$ be an extension of integral domains. Then $R=A+XB[X]$ is a particular 
  graded domain with the natural graduation. We reobtain [4, Corollary 1.2]:
\bigskip

\noindent{\bf Corollary 1.10.}\quad {\sl $Cl(A+XB[X])=HCl(A+XB[X])$ if and only if $B$ is
integrally closed.} \medskip 

\noindent {\bf Proof.}\quad This follows from Corollary 1.8 and the fact that $A+XB[X]$ is 
almost normal if and only if $B$ is integrally closed.  $\diamondsuit$ \bigskip 

Now, it seems natural to ask whether the equivalence ``$Cl(R)=HCl(R)$ $\Leftrightarrow$ $R$ 
is almost normal"  always holds for a graded
domain $R$. By the proof of (iii)$\Rightarrow$ (i) of Lemma 1.6, the implication 
``$Cl(R)=HCl(R)$ $\Rightarrow$ $R$ is almost normal" is always true.
However, the converse fails to be true in general as the following example
shows.\bigskip

\noindent{\bf Example 1.11.}\quad Let $K$ be a field and let $X$,$Y$, and $Z$
be three indeterminates over $K$. Set $T= K[X,Y,Z]/(YZ+X-X^2)$. Then $T$ is an integral
domain  and $T= K[x,y,z]$, where $yz=x(x-1)$. Let $d$ be an integer and set
$T_d=K[x]y^d$ if $d>0$, $T_0=K[x]$, and $T_d=K[x]z^{-d}$ if $d<0$. Then
$T=\oplus_{n\in Z} T_d$ is a $Z$-graded integral domain (cf. [18, p. 13]). 
Now, let $R=\oplus_{n\in Z} R_d$ be the $Z$-graded subring of $T$ defined as follows: $R_d=T_d$
if $d\not=0$ and $R_0= K+x(x-1)K[x]$. Then $R$ is almost normal and $Cl(R)\not=HCl(R)$. \medskip

Let $S$ be the multiplicatively closed subset of nonzero 
homogeneous elements of $R$. Then $R_S=K(x)[y,z]$ is a $Z$-graded
integral domain with $(R_S)_d=K(x)y^d$ if $d>0$, $(R_S)_0=K(x)$, and
$(R_S)_d=K(x)z^{-d}$ if $d<0$. Now, let $\xi\in (R_S)_d$ be integral over $R$ with
$d>0$. Then $\xi^n+F_{n-1}\xi^{n-1}+\cdots+F_0=0$ for some $F_0,\ldots, F_{n-1}\in R$.
Since $\xi$ is homogeneous of degree $d$, we may assume that $F_i\in R_{(n-i)d}$ for
each $i$. It follows that $\xi=\varphi(x)y^d$ for some $\varphi(x)\in K(x)$ and,  for
each $i$, $F_i=f_i(x)y^{(n-i)d}$ for some $f_i(x)\in K[x]$. Hence 
$\varphi^n+f_{n-1}\varphi^{n-1}+\cdots+f_0=0$. Since $K[x]$ is integrally closed, then
$\varphi(x)\in K[x]$, and hence $\xi\in R_d$. The case $d<0$ is similar. 
Therefore, $R$ is almost normal.\par
 Now we exhibit an invertible ideal of $R$ which is not proportional to any invertible
homogeneous ideal of $R$. Let $a=x^2(x-1)$, $b=x-y$, $c=x(x-1)-xy$, $d=1+(z/x)$,
$e=z+x-1$, and $f=x-1$. Set $I=(a,b,c)$ and $J=(d,e,f)$. Then $I$ and $J$ are two
fractional ideals of $R$ with $IJ\subset R$ and $1=
16af-(4x(x-1)-1)[(be-cd)^2-(bd)^2+2bd]\in IJ$. It follows that $I$ is an invertible
ideal of $R$. Now, assume that there exist $q\in R_S$ and $H$ a homogeneous
(integral) ideal of $R$ such that $I=qH$. Since $a\in I$ and $a$ is a  homogeneous
element of $R$, necessarily $q$ is homogeneous in $R_S$. Therefore, $I$ is a 
homogeneous fractional ideal of $R$. It follows that $I=C(I)=(x,y)$ and $J=C(J)=(1,z/x)$. Thus,
$IJ=(x,x-1,y,z)$, a contradiction since $x\notin R$. $\diamondsuit$ 
\bigskip

\noindent {\bf Example 1.12.}\quad Let $T=C[X]$, where $C$ is the field of 
complex numbers. Then $R=Z[\sqrt{-5}]+XC[X]$, graded in the natural way, 
is an almost normal graded domain which is not integrally closed (and hence 
a non Krull domain) with $R_0\subset R$ inert. In this case, 
$Cl(R)=HCl(R)= Cl(Z[\sqrt{-5}])=Z/2Z$ (cf. [5, Theorem 3.12(2)]). $\diamondsuit$
\vskip1cm

\noindent {\bf 2. Application: a decomposition theorem for semigroup rings }
\medskip

Let $A$ denote an integral domain with quotient field
K and $\Gamma$  a nonzero torsionless grading monoid with quotient group $G=<\Gamma>$. 
In the first part of this section, we focus on the specific case of semigroup rings. 
These constitute maybe the most important 
class of graded domains. We appeal to the main theorem of Section 1 to establish a decomposition theorem for the
class group of a semigroup ring. Specifically, we show that $Cl(A[\Gamma]) \cong Cl(A)\oplus 
HCl(K[\Gamma])$ if and only if $A[\Gamma]$ is integrally closed. This recovers most of the well-known results stated 
for the classic contexts of polynomial rings [11] and Krull semigroup rings [3].
The problem breaks into two parts: first, we explore the question of when $Cl(A[\Gamma])$ coincides with
$HCl(A[\Gamma])$; then state that  $HCl(A[\Gamma]) \cong Cl(A)\oplus HCl(K[\Gamma])$.\par \bigskip

 \noindent{\bf Lemma 2.1.}\quad {\sl Let $A[\Gamma]$ be any  semigroup ring. The following
statements are equivalent.\par 
(i) $A[\Gamma]$  is almost normal ;\par 
(ii) $A[\Gamma]$  is integrally closed ;\par 
(iii) $A$ and $\Gamma$ are integrally closed.}
\medskip 

\noindent {\bf Proof.}\quad For (i)$\Rightarrow$(ii), see the proof of  [1, corollary 3.9]. The implication 
(ii)$\Rightarrow$(i) is obvious. Finally for (ii)$\Leftrightarrow$(iii), see  [13, Corollary
12.11]. $\diamondsuit$  
\bigskip 

It is known that the extension $A\subset A[\Gamma]$ is always inert.
Thus, as a consequence of Theorem 1.1 and Lemma 2.1, we have:\bigskip 

\noindent{\bf Proposition 2.2.}\quad  {\sl Let $A[\Gamma]$ be a semigroup ring. Then
$Cl(A[\Gamma])=HCl(A[\Gamma])$ if and only if $A[\Gamma]$ is integrally closed. $\diamondsuit$}\bigskip

In order to prove the second (and remaining) part of our main theorem, Theorem 2.7, we need some 
preliminary results. We begin our discussion by handling a number of technical points. 
Let  ${\cal F}(\Gamma)$ denote the set of all fractional ideals of $\Gamma$. Under
ordinary addition of subsets of $G$, that is, $X+Y=\{x+y\mid x\in X \,\hbox{and}\,
y\in Y\}$, ${\cal F}(\Gamma)$ is a commutative monoid with zero element $\Gamma$. If
$Y, Z\in {\cal F}(\Gamma)$, then $(Y:Z)$ is defined to be the fractional ideal
$(Y:Z)=\{g\in G\mid\,g+Z\subset Y\}$. The fractional ideal $Y^{-1}=(\Gamma :Y)$ (resp.,
$Y_v=(Y^{-1})^{-1}$) is called the inverse (resp., the $v$-closure) of $Y$. We say that $Y$ is  
divisorial or  a $v$-ideal if $Y_v=Y$. The ideal $Y$ is  $v$-finite if 
$Y=(F+\Gamma)_v$ for some finite subset $F$ of $G$. Note that finitely generated
fractional ideals of $\Gamma$ are of the form $F+\Gamma$, where $F$ is a finite
subset of $G$.  For more details about  the $v$-operation on semigroups,  see [13, p. 215].\par
 
 Now, let $Y$ be a fractional ideal of $\Gamma$. We define
$Y_t=\cup\{(F+\Gamma)_v\mid\, F\subset Y$ a finite subset$\}$. It is easily seen that if
$\alpha\in G$ and $Y,Z\in {\cal F}(\Gamma)$, then:\medskip

(i) $(\alpha+\Gamma)_t=\alpha+\Gamma$; $(\alpha+Y)_t=\alpha+Y_t$.\par 
(ii) $Y\subset Y_t$; if $Y\subset Z$, then $Y_t\subset Z_t$.\par 
(iii) $(Y_t)_t=Y_t$.\medskip 

\noindent Therefore, $Y\mapsto Y_t$ defines  a $\ast$-operation on $\Gamma$ (cf. [19, section 10]),  
called  the $t$-operation. It is  done in analogy with the $t$-operation for domains. 
A fractional ideal $Y$  of $\Gamma$ is a $t$-ideal if
$Y_t=Y$, or equivalently for each finite subset $F\subset Y$, $(F+\Gamma)_v\subset Y$.
Clearly, if $Y$ is a fractional ideal of $\Gamma$, $Y_t\subset Y_v$ and hence a $v$-ideal is a $t$-ideal. 
Let $t(\Gamma)$  denote the subset of ${\cal F}(\Gamma)$ of all $t$-ideals of $\Gamma$. One may 
easily see that $(Y+Z)_t=(Y_t+Z)_t=(Y_t+Z_t)_t$ for 
$Y,Z$ fractional ideals of $\Gamma$, and hence $t(\Gamma)$ forms a commutative
monoid,  with zero element $\Gamma$, under the operation
$(Y,Z)\mapsto(Y+Z)_t$.  An invertible element for this operation  is called a
$t$-invertible (fractional) $t$-ideal of $\Gamma$. As in the case of domains, 
a $t$-invertible  $t$-ideal is always divisorial and $v$-invertible. Conversely, 
a divisorial $v$-invertible  ideal $Y$ of $\Gamma$ is  $t$-invertible   if and only if  
$Y$ and $Y^{-1}$ are $v$-finite.\par

 In what follows, the graduation on the semigroup ring $A[\Gamma]$ will always mean the
natural graduation.\bigskip 

\noindent{\bf Lemma 2.3.}\quad {\sl Let $A[\Gamma]$ be a semigroup ring and let $I,J$
(resp., $Y,Z$ ) be two fractional ideals of $A$ (resp., $\Gamma$). Then\par

 (1) $I[Y]$ is a homogeneous fractional ideal of $A[\Gamma]$.\par 
(2) $I[Y]$ is finitely generated if and only if $I$ and $Y$ are finitely
generated.\par   
If, moreover, $I$ and  $J$ are nonzero, then\par 
(3) $(I[Y]:J[Z])=(I:J)[(Y:Z)]$.\par 
(4) $(I[Y])_v=I_v[Y_v]$.\par 
(5) $(I[Y])_t=I_t[Y_t]$. }  \medskip 

\noindent {\bf Proof.}\quad (1) First, note that $I[Y]$
is the subset of elements of $K[G]$ of the form $\sum a_iX^{\alpha_i}$, where $a_i\in
I$ and $\alpha_i\in Y$. We have $A[\Gamma]I[Y]\subset I[Y+\Gamma]= I[Y]$. On the other hand, if
$0\not=a\in A$ and $\alpha \in \Gamma$ are such that $aI\subset A$ and
$\alpha+Y\subset\Gamma$, then $aX^{\alpha}I[Y]\subset A[\Gamma]$. Hence $I[Y]$ is a
fractional ideal of $A[\Gamma]$. The fact that $I[Y]$ is homogeneous  is obvious.\par
(2)   If $I$ and $Y$ are finitely generated ideals, it is obvious that $I[Y]$ is
finitely generated. Conversely, suppose that $I[Y]$ is finitely generated. Since
$I[Y]$ is homogeneous, there exist $a_1,\ldots, a_n\in I$ and $\alpha_1, \ldots,
\alpha_n\in Y$ such that $I[Y]=(a_1X^{\alpha_1}, \ldots, a_nX^{\alpha_n})$. Now, one may
easily check that $I=(a_1,\ldots, a_n)$ and $Y=\{\alpha_1, \ldots,
\alpha_n\}+\Gamma$.\par 

(3) Clearly, $(I:J)[(Y:Z)]\subset (I[Y]:J[Z])$. For the reverse inclusion, let $f\in(I[Y]:J[Z])$. 
So $fJ[Z]\subset I[Y]$ and then $f\in (I:J)[G]$. On
the other hand, $fJ[Z]\subset I[Y]$ implies that $\alpha +Z\subset Y$ for each
$\alpha\in Supp(f)$. Hence $Supp(f)\subset(Y:Z)$. It follows that
$f\in(I:J)[(Y:Z)]$.\par 
(4) Follows from (3).\par 
(5) By (2) and (4), we have
$$\eqalign{(I[Y])_t&=\bigcup\{(F[T])_v\mid\, F\subset I \,\hbox{and}\, T\subset Y,
F\,\hbox{and}\,T\, \hbox{are ideals of finite type}\}\cr &=\bigcup\{F_v[T_v]\mid\,
F\subset I \,\hbox{and}\, T\subset Y, F\,\hbox{and}\,T\,\hbox{are  of finite
type}\}\cr &=I_t[Y_t]. \diamondsuit \cr }$$\bigskip 

\noindent{\bf Corollary 2.4.}\quad {\sl Let $A[\Gamma]$ be a semigroup ring, $I$ a nonzero
fractional ideal of A, and $Y$ a fractional ideal of $\Gamma$. Then\par (1) $I[Y]$ is a
$v$-ideal (resp., $t$-ideal) if and only if $I$ and $Y$ are $v$-ideals (resp.,
$t$-ideals).\par (2) $I[Y]$ is $v$-finite if and only if $I$ and $Y$ are $v$-finite.\par
(3) $I[Y]$ is $v$-invertible (resp., $t$-invertible) if and only if
$I$ and $Y$ are $v$-invertible (resp., $t$-invertible).}\medskip

\noindent {\bf Proof.}\quad (1) is a
consequence of the statements (4) and (5) of Lemma 2.3.\par
(2) Follows from  (2) and (4) of Lemma 2.3.\par  
(3) Follows from the fact that
$(I[Y](I[Y])^{-1})_v=(II^{-1})_v[(Y+Y^{-1})_v]$ and \par\noindent
$(I[Y](I[Y])^{-1})_t\,\,=(II^{-1})_t[(Y+Y^{-1})_t]$.  $\diamondsuit$ \bigskip 

Next we give a characterization of the homogeneous divisorial ideals of a semigroup
ring.\bigskip 

\noindent{\bf Proposition 2.5.}\quad {\sl Let $A[\Gamma]$ be a semigroup ring. The following statements are equivalent.\par 
(i) $I$ is a homogeneous fractional $v$-ideal of $A[\Gamma]$ ;\par
 (ii) $I=J[Y]$ for some fractional $v$-ideals $J$ and $Y$ of $A$ and $\Gamma$, respectively.}\medskip

\noindent {\bf Proof.}\quad (i)$\Rightarrow$(ii) Let $I$ be
a nonzero homogeneous fractional ideal of $A[\Gamma]$. Then there exist $0\not=c\in A$ and $\gamma \in \Gamma$
such that $cX^{\gamma}I\subset A[\Gamma]$, so $I\subset K[G]$. Let $Y$ be the set of
degrees of all homogeneous elements of $I$ and let $J$ be  the $A$-submodule of $K$
generated by the coefficients of all elements of $I$. We have $Y+\Gamma\subset Y$,
$\gamma+Y\subset\Gamma$ and $cJ\subset A$. Hence $Y$ and $J$ are fractional ideals.
Next we show that $I=J[Y]$. The inclusion $I\subset J[Y]$ is trivial. For the reverse 
inclusion, let $f, 0\not=g\in A[\Gamma]$ such that $I\subset {f\over g}A[\Gamma]$. Let 
$aX^{\alpha}\in I$ with $a\not=0$. Then ${f\over g}={{aX^{\alpha}}\over h}$ for 
some $0\not=h\in A[\Gamma]$. Now, let $bX^{\beta}\in I$. Then 
$bX^{\beta}\in {{aX^{\alpha}}\over h}A[\Gamma]$. That is, 
$bX^{\beta}h\in aX^{\alpha}A[\Gamma]$, so $bh\in aA[\Gamma]$.
Therefore,  $bX^{\alpha}\in {f\over g}A[\Gamma]$ and hence $bX^{\alpha}\in I$ 
(since $I$ is divisorial). Hence $J[Y]\subset I$ and $I=J[Y]$.\par
(ii)$\Rightarrow$(i) Follows from Corollary 2.4.  $\diamondsuit$ \bigskip

\noindent{\bf Theorem 2.6.}\quad {\sl Let $A[\Gamma]$ be a semigroup ring. We have the
following splitting exact sequence of natural homomorphisms:} $$0\rightarrow
Cl(A)\rightarrow HCl(A[\Gamma])\rightarrow HCl(K[\Gamma])\rightarrow 0$$\medskip

\noindent {\bf Proof.}\quad Since $A[\Gamma]$ is a flat $A$-module, the natural homomorphism 
$Cl(A){\!}\rightarrow{\!} Cl(A{\!}[\Gamma])$,  $[J]\mapsto[J[\Gamma]]$ is  well-defined 
(cf. [5, Proposition 2.2]), and it induces a natural homomorphism 
$\phi: Cl(A){\!}\rightarrow{\!} HCl(A{\!}[\Gamma])$. On the other hand, since 
$K[\Gamma]$ is a  quotient ring of $A[\Gamma]$, we have the natural 
homomorphism  $Cl(A[\Gamma])\rightarrow
Cl(K[\Gamma])$, $[I]\mapsto [IK[\Gamma]]$. It induces
a natural homomorphism \hbox{$\psi:\, HCl(A[\Gamma])\rightarrow HCl(K[\Gamma])$}.
 By Corollary 2.4 and Proposition 2.5, if  $I$ is a
homogeneous  $t$-invertible $t$-ideal of $A[\Gamma]$,  $I=J[Y]$ for some $t$-invertible fractional 
$t$-ideals $J$ and $Y$ of $A$ and $\Gamma$, respectively. Hence
$\psi([I])=\psi([J[Y]])=[K[Y]]$. Thus we have the sequence:  $${Cl(A)\buildrel
\phi \over {\rightarrow}} {HCl(A[\Gamma])\buildrel \psi \over {\rightarrow} HCl(K[\Gamma])}$$
 
 \indent We first show that $\phi$ is injective. Let $J$ be a $t$-invertible $t$-ideal of 
 $A$ such that $J[\Gamma]=uA[\Gamma]$ for some homogeneous element $u\in K[G]$. Then 
 $u\in A$ and $J=uA$. Hence $\phi$ is injective.\par
Next we show that $Im\phi= Ker\psi$. Clearly $Im\phi\subset Ker\psi$. To show that 
$Ker\psi\subset Im\phi$, let $I$ be a homogeneous $t$-invertible $t$-ideal of
$A[\Gamma]$ such that $IK[\Gamma]=fK[\Gamma]$ for some homogeneous element $f\in K[\Gamma]$. 
We may assume that $f=X^{\alpha}$ for some $\alpha\in \Gamma$. Now, set
$I_1=X^{-\alpha}I$. Then $I_1\subset A[\Gamma]$ and $I_1K[\Gamma]=K[\Gamma]$. It
follows that $J=I_1\cap A\not=0$ and by [1, Proposition 5.7], $I_1=J[\Gamma]$. Since
$I_1$ is a $t$-invertible $t$-ideal, $J$ is a $t$-invertible $t$-ideal of $A$ (cf.
Corollary 2.4). Hence $\phi([J])=[I]$ and $Ker\psi\subset Im\phi$.\par
 Now, let $I$ be a homogeneous $t$-invertible $t$-ideal of $K[\Gamma]$, by Proposition 2.5,
$I=K[Y]$ for some ideal $Y$ of $\Gamma$. From Corollary 2.4, $Y$ is a $t$-invertible
$t$-ideal of $\Gamma$ and $A[Y]$ is a $t$-invertible $t$-ideal of $A[\Gamma]$. Now, consider the map $$\psi': HCl(K[\Gamma])\rightarrow HCl(A[\Gamma]),\, [I]=[K[Y]]\mapsto
[A[Y]].$$ It is clear that $\psi'$ is a well-defined homomorphism and $\psi o\psi'=i$,
the identity map. $\diamondsuit$ \bigskip

Finally, we are able to announce our decomposition theorem.\bigskip  

\noindent{\bf Theorem 2.7.}\quad {\sl Let  $A[\Gamma]$ be a semigroup ring. Then 
$Cl(A[\Gamma])\cong Cl(A)\oplus HCl(K[\Gamma])$ if and only if $A[\Gamma]$ is
integrally closed.} \medskip

\noindent {\bf Proof.}\quad  It follows from Proposition 2.2 (i.e., Theorem 1.1) and Theorem 2.6.	 $\diamondsuit$ 
\bigskip

The following corollary is a straightforward consequence of the above results.\bigskip

\noindent{\bf Corollary 2.8.} \quad {\sl Let $A[\Gamma]$ be an integrally closed semigroup ring. Then 
$Cl(A[\Gamma])\cong Cl(A)\oplus Cl(K[\Gamma])$.  $\diamondsuit$} 
\bigskip

The second part of this section is devoted to semigroups. Here our ambiant 
semigroup $\Gamma$ is actually a nonzero torsionless grading monoid. This is indispensable to ensure that the 
associated semigroup ring $A[\Gamma]$ will be an integral domain. Our purpose is to extend Chouinard's 
results on Krull semigroups (cf. [8]) to arbitrary semigroups. As an application, we establish
an interesting result, Theorem 2.12, on the natural homomorphism 
$\phi: Cl(A)\rightarrow Cl(A[\Gamma])$ along with a few consequences. \par

In [8, Lemma 1, p.1463], the author proved  that 
$Cl(K[\Gamma])\cong Cl(\Gamma)$ for any field $K$ and any Krull semigroup $\Gamma$
with $\Gamma\cap(-\Gamma)=0$. This agreed with the fact that $Cl(K[\Gamma])$ is
not controlled by $K$ provided $K[\Gamma]$ is a Krull domain [3, Proposition 7.3(2)]. \par

In analogy with 
the case of integral domains, we define the  ($t$-)class group of the semigroup $\Gamma$, 
denoted $Cl(\Gamma)$, to be the group of $t$-invertible fractional
$t$-ideals of $\Gamma$ under $t$-multiplication modulo its subgroup of principal fractional  ideals. 
Also, in a  Krull semigroup, the $t$-operation and the $v$-operation coincide. 
Thus, if $\Gamma$ is a Krull semigroup, $Cl(\Gamma)$ is just the divisor class group 
of $\Gamma$ defined in [8], see also [13, Section 16].\bigskip

 \noindent{\bf Theorem 2.9.}
\quad {\sl Let $K$ be a field. Then $HCl(K[\Gamma])\cong
Cl(\Gamma)$.}\medskip

\noindent {\bf Proof.}\quad Consider the map $\varphi:\, Cl(\Gamma)\rightarrow
HCl(K[\Gamma]),\, [Y]\mapsto [K[Y]].$ If $Y$ is a $t$-invertible $t$-ideal of $\Gamma$,
then by Corollary 2.4, $K[Y]$ is a homogeneous $t$-invertible $t$-ideal of
$K[\Gamma]$. Hence $\varphi$ is a well-defined homomorphism. Now, let $Y$ be a
$t$-invertible $t$-ideal of $\Gamma$ such that $K[Y]=fK[\Gamma]$ for some homogeneous 
element $f\in K[G]$. We may suppose that $f=X^{\alpha}$ for
some $\alpha\in Y$. That is, $K[Y]=X^{\alpha}K[\Gamma]$. Therefore
$K[Y]=K[\alpha+\Gamma]$, and hence $Y=\alpha+\Gamma$ is principal. It follows that $\varphi$
is injective. To show that $\varphi$ is also surjective, let $I$ be a homogeneous
$t$-invertible $t$-ideal of $K[\Gamma]$. By Proposition 2.5 and Corollary 2.4, there
exists $Y$, a $t$-invertible $t$-ideal of $\Gamma$, such that $I=K[Y]$. Thus $\varphi([Y])= [I]$. $\diamondsuit$ \bigskip 

As a consequence of Theorem 2.9, we have the following corollaries which recover  [8, Lemma 1].
\bigskip

\noindent{\bf Corollary 2.10.} \quad {\sl Let $K$ be a field. 
If $\Gamma$ is integrally closed, then $Cl(K[\Gamma])\cong
Cl(\Gamma)$.} \medskip 

\noindent {\bf Proof.}\quad It follows from Lemma 2.1, Proposition 2.2, and Theorem 2.9. $\diamondsuit$ \bigskip
 
\noindent{\bf Corollary 2.11.} \quad {\sl Let $A[\Gamma]$ be an integrally closed semigroup ring. Then 
$Cl(A[\Gamma])\cong Cl(A)\oplus Cl(\Gamma)$}.
\medskip

\noindent {\bf Proof.}\quad It follows from Theorem 2.7 and Theorem 2.9. $\diamondsuit$ \bigskip

We close this section by a brief study of the canonical homomorphism 
$\phi:  Cl(A)\rightarrow Cl(A[\Gamma])$, $[J]\mapsto [J[\Gamma]]$. Let $R$ be a graded 
domain. It is well-known  that, when defined, the 
homomorphism $\phi: Cl(R_0)\rightarrow Cl(R)$, $[I]\mapsto [(IR)_t]$ is not an isomorphism in
general. For example, see [3, section 6]. Nevertheless, in the case of a semigroup 
ring we give a complete characterization for $\phi$ to be an isomorphism.\bigskip
  
 \noindent{\bf Theorem 2.12.}\quad
{\sl Let $A$ be an integral domain. Then $Cl(A[\Gamma])=Cl(A)$ if
and only if $A$ and $\Gamma$ are integrally closed and $Cl(\Gamma)=0$.} \medskip

\noindent {\bf Proof.}\quad Assume that $Cl(A[\Gamma])=Cl(A)$. Since $\phi$ maps into $HCl(A[\Gamma])$,
then $Cl(A[\Gamma])=HCl(A[\Gamma])$. Therefore,  $A[\Gamma]$ is integrally closed (cf.
Proposition 2.2) and   $Cl(\Gamma)=0$ (cf. Corollary 2.11). The converse follows from Lemma 2.1 and  
Corollary 2.11. $\diamondsuit$ \bigskip

We conclude with some corollaries and examples illustrating (the scope of) Theorem 2.12.
 \bigskip
  
Theorem 2.12 generalizes [11, Theorem 3.6] (see also Corollary 1.9). To see this,
let $\Gamma= \oplus Z_+e_\alpha$. Clearly, 
$G=\oplus Ze_\alpha$. Since $\Gamma$ is factorial (cf. [13, Theorem 6.8]),
then $Cl(\Gamma)= 0$. Hence 
  \bigskip

 \noindent{\bf Corollary 2.13.}\quad {\sl $Cl(A[\{X_\alpha\}] = Cl(A)$ if and only if $A$ is
integrally closed. $\diamondsuit$} \bigskip
 
 As an other consequence of
Theorem 2.12, we have the following result on group rings which recovers [16, Proposition 5.3].\bigskip

 \noindent{\bf Corollary 2.14.}\quad {\sl Let $A$ be
an integral domain and $G$ a nonzero torsion-free abelian group. Then $Cl(A[G])=Cl(A)$ if and
only if $A$ is integrally closed. $\diamondsuit$} 
\bigskip

 \noindent{\bf Corollary 2.15.}\quad {\sl Let $A$ be
an integral domain. Then 
$Cl(A[\{X_\alpha, X_\alpha^{-1}\}])=Cl(A)$ if and
only if $A$ is integrally closed.  $\diamondsuit$} 
\bigskip

\noindent{\bf Example 2.16.}\quad Let $\Gamma=Z\times Z_+$ and let $A$ be an integrally 
 closed domain. By [13, Theorem 6.8], $\Gamma$ is factorial. Hence $\Gamma$ is 
 integrally closed and $Cl(\Gamma)=0$. Thus $Cl(A[\Gamma])=Cl(A[X, X^{-1}, Y])=Cl(A)$.	 $\diamondsuit$ 
\bigskip
 
 \noindent{\bf Example 2.17.} \quad Let $\Gamma=\cup_{n\ge0}{1\over{p^n}}Z_+$, where $p$ is a positive prime 
integer. Since $\Gamma$ does not satisfy 
the a.c.c on subsemigroups, then $\Gamma$ is not a Krull semigroup. Now let 
$A=Z[\sqrt{-5}]$ and $K$ its quotient field. By [13, Theorem 13.5], $K[\Gamma]$ 
is a Bezout domain. Hence  $\Gamma$ is integrally closed and $Cl(\Gamma)=Cl(K[\Gamma])=0$. 
Thus $Cl(Z[\sqrt{-5}][\Gamma])=Cl(Z[\sqrt{-5}])=Z/2Z$.  $\diamondsuit$ 
\bigskip

Theorem 2.12 shows that the class group of $A[\Gamma]$ measures the failure of integral 
closure for the ring $A$ and for the semigroup $\Gamma$. The following example illustrates
this fact. It also illustrates the failure of Corollary 2.10 for non-integrally closed semigroups. 
\bigskip

\noindent{\bf Example 2.18.}\quad Let $\Gamma=\{0, 2, 3, 4,\ldots\}$ and $K$ be any field. 
Then $K[\Gamma]=K[X^2, X^3]$ and $\Gamma$ is not integrally closed. Let 
$Y_0=\{2, 3, 4, \ldots\}$. Then  $Y_0$ is an integral ideal of $\Gamma$, 
and if $Y$ is a nonprincipal integral ideal of $\Gamma$,  $Y=n+Y_0$, for some integer
 $n\ge 0$. On the other hand, one can easily  check that $Y_0$ is a divisorial 
 ideal of $\Gamma$. Hence all ideals of $\Gamma$ are divisorial. Thus 
 $Cl(\Gamma)=Pic(\Gamma)$. Now let $Y$ be an invertible ideal of $\Gamma$. 
 Then $Y+(\Gamma:Y)=\Gamma$. Let $n\in Y$ and $m\in (\Gamma:Y)$ such that $n+m=0$. 
 We have $Y=n+m+Y\subset n+\Gamma$, so $Y=n+\Gamma$ is a principal ideal. 
 Hence $Cl(\Gamma)=Pic(\Gamma)=0$. Thus $Cl(\Gamma)=HCl(K[\Gamma])=0$, 
 while $Cl(K[\Gamma])=K$ (cf. [5, Example 4.7(2)]). $\diamondsuit$ 
\vskip 1cm

{\bf References} \medskip

\item{[1]} D.D. Anderson, D.F. Anderson, Divisorial ideals and invertible ideals in a 
graded integral domain, J. Algebra 76 (2) (1982) 549-569.\medskip

\item{[2]} D.D. Anderson, D.F. Anderson, Divisibility properties of graded 
domains, Canad. J. Math.  34 (1) (1982) 196-215.\medskip

\item{[3]} D.F. Anderson, Graded Krull domains, Comm. Algebra 7 (1979) 79-106.\medskip

\item{[4]} D.F. Anderson, S. El Baghdadi, S. Kabbaj,
The homogeneous class group of A+XB[X] domains, International Journal of Commutative Rings (to appear). \medskip

\item{[5]} D.F. Anderson, A. Ryckaert, The class group of $D+M$, 
J. Pure Appl. Algebra 52 (1988) 199-212.\medskip

\item{[6]} A. Bouvier, Le groupe des classe d'un anneau int\`egre, 
in: 107\` eme congr\`es des soci\'et\'es savantes, vol. IV, Brest, 1982, pp. 85-92.\medskip

\item{[7]} H. Bass, M. P. Murthy, Grothendieck groups and Picard groups 
of abelian group rings, Ann. of Math. 86 (1967) 16-73.\medskip

\item{[8]} L.G. Chouinard, Krull semigroups and divisors class 
groups, Canad. J. Math.  33 (6) (1981) 1459-1464.\medskip

\item{[9]} P. M. Cohn, Bezout rings and their subrings,
Proc. Cambridge Philos. Soc. 64 (1968) 251-264.\medskip 

\item{[10]} R.M. Fossum, The divisor class group of a Krull domain, 
Springer-Verlag, New York, 1973.\medskip

\item{[11]} S. Gabelli, On divisorial ideals in polynomial rings over 
Mori domains, Comm. Algebra 15  (1987) 2349-2370.\medskip

\item{[12]} R. Gilmer, Multiplicative ideal theory, Marcel 
Dekker, New York, 1972.\medskip

\item{[13]} R. Gilmer, Commutative semigroup rings, Chicago Lecture 
Notes in Mathematics, Chicago, 1984.\medskip

\item{[14]} R. Gilmer, R. Heitmann, On Pic(R[X]) for R seminormal, 
J. Pure Appl. Algebra 16  (1980) 251-257. 
\medskip

\item{[15]} E.G. Houston, M. Zafrullah, On $t$-invertibility II, 
Comm. Algebra 17 (8) (1989) 1955-1969. \medskip

\item{[16]} R. Matsuda, On algebraic properties of infinite group rings, 
Bull. Fac. Sci. Ibaraki Univ., Series A. Math., 7  (1975) 29-37.\medskip

\item{[17]} R. Matsuda, On the content condition of a graded integral domain, 
Comment. Math. Univ. St. Paul. 33 (1) (1984) 79-86.\medskip

\item{[18]} R. Matsuda, $\phi: Pic(R_0)\rightarrow HPic(R)$ is not necessarily surjective, 
Bull. Fac. Sci. Ibaraki Univ., Series A. Math., 16 (1984) 13-15.\medskip

 \item{[19]} R. Matsuda, Torsion-free abelian semigroup rings VI, 
 Bull. Fac. Sci. Ibaraki Univ., Series A. Math., 18 (1986) 23-43.\medskip

 \item{[20]} D.G. Northcott, A generalization of a theorem on the content of
polynomials, Proc. Cambridge Philos. Soc. 55 (1959) 282-288.\medskip 

\item{[21]} D.G. Northcott, Lessons on rings, modules, and multiplicities, Cambridge
Univ. Press, Cambridge, 1968.\medskip

\item{[22]} C. Pedrini, On the $K_0$ of certain polynomial extensions, 
Lecture Notes in Mathematics, vol. 342, Springer-Verlag, New York, 1973, pp. 92-108.\medskip
\end